\documentclass[12pt,reqno]{amsart}

\numberwithin{equation}{section}

\usepackage[final, 
color]{showkeys}
\definecolor{refkey}{gray}{.85}
\definecolor{labelkey}{gray}{.85}

\usepackage{AKstyle}

\usepackage{caption}
\usepackage[labelfont=rm]{subcaption}

\usepackage[pagebackref=true, colorlinks]{hyperref}

\hypersetup{pdffitwindow=true,linkcolor=blue,citecolor=blue,urlcolor=blue,menucolor=blue}

\usepackage{comment}

\makeatletter
\let\orgdescriptionlabel\descriptionlabel
\renewcommand*{\descriptionlabel}[1]{%
  \let\orglabel\label
  \let\label\@gobble
  \phantomsection
  \edef\@currentlabel{#1}%
  \let\label\orglabel
  \orgdescriptionlabel{#1}%
}
\makeatother

\begin{document}

\author{Jean Bourgain}
\thanks{JB is partially supported by NSF grant DMS-1301619.}
\email{bourgain@ias.edu}
\address{School of Mathematics, IAS, Princeton, NJ}

\author{Alex Kontorovich}
\thanks{AK is partially supported by
an NSF CAREER grant DMS-1254788 and  DMS-1455705, an NSF FRG grant DMS-1463940, an Alfred P. Sloan Research Fellowship, and a BSF grant.}
\email{alex.kontorovich@rutgers.edu}
\address{Department of Mathematics, Rutgers University, New Brunswick, NJ}

\title[Reciprocal Geodesics]
{Beyond Expansion III: 
Reciprocal
Geodesics}

\begin{abstract}
We prove the existence of infinitely many low-lying 
and
fundamental 
closed geodesics on the modular surface which are
reciprocal, 
that is,
invariant
under time reversal. 
The method combines ideas from Parts I and II of this series, namely the dispersion method in bilinear forms, as applied to thin semigroups coming from 
restricted
continued fractions.
\end{abstract}
\date{\today}
\subjclass[2010]{}
\maketitle
\tableofcontents


\section{Introduction}\label{sec:intro}

We quote from Sarnak's lecture  \cite{Sarnak2010}  
regarding the genesis of the Affine Sieve \cite{BourgainGamburdSarnak2006, BourgainGamburdSarnak2010, SalehiSarnak2013}:

\vskip.1in
{\small ``For me the starting point of this investigation was in 2005 when Michel and Venkatesh asked me about the existence of poorly distributed closed geodesics on the modular surface. It was clear that Markov's constructions of his geodesics using his Markov equation provided what they wanted but they preferred quadratic forms with square free discriminant. This raised the question of sieving in this context of an orbit of a group of (nonlinear) morphisms of affine space.''}

The initial question arose in Einsiedler-Lindenstrauss-Michel-Venkatesh's investigations into higher rank analogues  of Duke's theorem \cite{Duke1988}, and asked (see the discussion below \cite[Thm 1.10]{ELMV2009}) for an infinitude of low-lying (that is, being poorly distributed by not entering the cusp) {\it fundamental} geodesics (i.e., those corresponding to fundamental classes of binary quadratic forms).
This problem was
 solved in Part II of our series; see \cite{BourgainKontorovich2015, Kontorovich2016} for a detailed  discussion. But the question of 
an infinitude of fundamental Markov 
geodesics (for a discussion of Markov geodesics, see, e.g., \cite[p. 226]{Sarnak2007}) remains wide open, despite recent progress on the ``strong approximation'' aspect  in \cite{BourgainGamburdSarnak2015, BourgainGamburdSarnak2016}. Such geodesics 
are all reciprocal, that is, equivalent to themselves under time-reversal of the geodesic flow. In this paper we relax Markov geodesics to just low-lying ones, and solve the problem of producing an infinitude of low-lying, funamental, reciprocal geodesics.

\subsection{Statement of the Main Theorem}\

Before stating our main result, we give precise definitions of low-lying, fundamental, and reciprocal. By closed geodesic, we always mean primitive.

\begin{Def}
Given a compact subset $\cY$ of the unit tangent bundle of the modular surface 
$$
\cX=T^{1}(\PSL_{2}(\Z)\bk\bH)\cong\PSL_{2}(\Z)\bk\PSL_{2}(\R), 
$$
a closed geodesic $\g$ on $\cX$ is called {\bf low-lying} (with respect to $\cY$) if 
$
\g\subset\cY.
$
\end{Def}

\begin{Def}
As is well-known, closed geodesics on $\cX$ are in 1-1 correspondence with primitive conjugacy classes of hyperbolic elements of $\PSL_{2}(\Z)$, and also with equivalence classes of indefinite binary quadratic forms (see, e.g., \cite{Kontorovich2016}). The latter 
have 
discriminants, and we say that a closed geodesic has discriminant $D$ if its corresponding class does. 
The {\bf trace} of a closed geodesic is that of its corresponding conjugacy class.
Recall that a non-square discriminant $D$ is called fundamental if it is the discriminant of the real quadratic field $\Q(\sqrt D)$. We call a closed geodesic {\bf fundamental} if its discriminant is.
\end{Def}

\begin{Def}
The time reversal symmetry on $\cX$ corresponds to replacing all tangent vectors by their negatives; if a closed geodesic is invariant under this involution, it is called {\bf reciprocal}.
\end{Def}

Recall that the total number of all primitive closed geodesics, ordered by trace (which is equivalent to ordering by length), has the following well-known asymptotic: 
$$
\#\{\text{\it closed geodesics with trace}<X\}
\
\sim
\
{X^{2}\over 2\log X}.
$$
There are about square-root as many reciprocal geodesics, which makes intuitive sense, as the geodesic has to spend the second half of its life undoing the twists of its first half.

\begin{thm}[{Sarnak \cite[Thm. 2]{Sarnak2007}}]
$$
\#\{\text{reciprocal geodesics with trace}<X\}
\
\sim
\
\frac38 X.
$$
\end{thm}

Our main result produces almost as  many low-lying, fundamental, reciprocal geodesics.

\begin{thm}[Main Theorem]\label{thm:1}
For any $\eta>0$, there is a compact subset $\cY=\cY(\eta) \subset\cX$ so that 
$$
\#\{\text{low-lying, fundamental, reciprocal geodesics with trace}
<X\}
$$
$$
\hskip2in
\gg_{\eta}
\
 X^{1-\eta}.
$$
\end{thm}

\subsection{Ingredients}\

As in  Part II of our series \cite{BourgainKontorovich2015}, we must study restricted continued fractions, and to understand these, we use the semigroup 
\be\label{eq:GcADef}
\G_{\cA} \ := \
\<
\mattwo a110:a\le \cA
\>^{+}\ \cap \ \SL_{2}
,
\ee
of even length words in the generators displayed.
Write $B_{N}$ for the archimedean ball in $\SL_{2}(\R)$ with respect to the Frobenius metric:
$$
B_{N} \ := \
\{ g=\mattwos abcd \in \SL_{2}(\R)  \ : \ \tr(g^{\dag} g)=a^{2}+b^{2}+c^{2}+d^{2}<N^{2}\}.
$$
Hensley \cite{Hensley1989} estimates the size of an archimedean ball  in $\G_{\cA}$ to be
\be\label{eq:Hensley}
\#\G_{\cA}\cap B_{N} \ \asymp \ N^{2\gd_{\cA}},
\ee
where $\gd_{\cA}$ is the Hausdorff dimension of the limiting Cantor set,
$$
\fC_{\cA}\ : = \ 
\{
[0,a_{1},a_{2},\dots] \ : \ a_{j}\le \cA\text{ for all } j
\}.
$$
Here we are using the standard notation $x=[a_{0},a_{1},a_{2},\dots]$ for the continued fraction
$$
x=a_{0}+\cfrac{1}{a_{1}+\cfrac{1}{a_{2}+\ddots}}.
$$
These fractal dimensions are known to tend to $1$ as $\cA\to\infty$; indeed Hensley \cite{Hensley1992} has shown that:
\be\label{eq:gdcA}
\gd_{\cA} \ = \ 1- {6\over \pi^{2}\cA}+o\left(\frac1\cA\right).
\ee

The following lemmata give sufficient conditions for a closed geodesic -- represented by a hyperbolic conjugacy class $[\g]$ with $\g\in\SL_{2}(\Z)$ -- to be fundamental and reciprocal.
\begin{lem}[{\cite[Lemma 1.14]{BourgainKontorovich2015}}]\label{lem:fundGeo}
A sufficient condition for a closed geodesic $[\g]$ to be fundamental is that 
\be\label{eq:fundGeo}
\tr(\g)^{2}-4\text{ is square-free.}
\ee
\end{lem}

\begin{lem}[{See \cite{Sarnak2007}}]\label{lem:recipGeo}
A sufficient condition for a closed geodesic $[\g_{1}]$ to be reciprocal is that  it is
of the from  $\g_{1}=\g{}^{\dag}\g$, for some $\g\in\SL_{2}(\Z)$.
\end{lem}

We then reduce \thmref{thm:1} to the following sieving result.

\begin{thm}\label{thm:main}
For any $\eta>0$, there is an $\cA=\cA(\eta)<\infty$ so that:
$$
\#\{
\g\in\G_{\cA}\cap B_{N} \ : \
\tr(\g^{\dag}\g)^{2}-4\text{ is square-free}
\}
\ \gg \
N^{2 - \eta}
.
$$
\end{thm}

\begin{rmk}\label{rmk:AbsSpecGap}
As in Part II \cite{BourgainKontorovich2015}, we cannot simply execute the Affine Sieve, because the ``spectral gap'' is insufficiently robust relative to the growth exponent $\gd_{A}$, and we must produce an ``exponent of distribution'' going beyond that arising from expansion alone; see \rmkref{rmk:NeedAbsGap}.
To do this, we again create certain ``bilinear forms,'' and substitute ``resonance'' harmonics with abelian theory, which is much more tractable.
Unlike Part II, the direct approach fails
due to the nature of the quadratic forms arising in the error terms, and a version of Linnik's ``dispersion method'' is needed. Fortunately, such was just  developed in the 
``orbital sieve''
 context in Part I of our series \cite{BourgainKontorovich2015a}, and this comes to the rescue here.
\end{rmk}

\begin{rmk}
The main result in  Part II 
was proved unconditionally 
but would also follow immediately from a certain
``Local-Global Conjecture for thin orbits,''  see the discussion in \cite{Kontorovich2016}.
In contradistinction, \thmref{thm:1} does {\it not} follow from this conjecture,
because 
the function 
$$
\SL_{2}(\Z)\to \Z \ : \ \g\mapsto \tr(\g^{\dag}\g)
$$
 is quadratic in the entries, so cannot be onto when restricted to any $\G_{\cA}$; the image is itself thin! (For a definition of thinness in this context, see \cite[p. 954]{Kontorovich2014}.)
\end{rmk}

\subsection{Organization}\

The rest of the paper is organized as follows.  After some preliminary calculations in \secref{sec:prelim}, we state the Sieving Theorem and construct the bilinear forms in \secref{sec:setup} before analyzing the ``main term'' in
\secref{sec:main}. The error terms are analyzed in \secref{sec:err}, after which the Sieving Theorem is proved in \secref{sec:pfThm1}. Finally, putting together the above ingredients, we prove \thmref{thm:1} in \secref{sec:pfMain}.

\subsection{Notation}\

The transpose of a matrix $\g$ is written $^{\dag}\g$.
When a calculation involves modular arithmetic, an overbar, $\bar a$, shall denote the multiplicative inverse of $a$.
The constants $C, c$ are absolute but may change from line to line.
We use the notation $f\ll g$ and $f=O(g)$ to mean $f(x)\le C g(x)$ for all $x>C$, where $C$ is an  implied constant. We write $f\asymp g$ for $g\ll f \ll g$.
Unless otherwise specified, implied constants depend at most on $\cA$, which is treated as fixed, and possibly on an arbitrarily small $\vep>0$.

\subsection*{Acknowledgments}\

It is our pleasure to thank
 Nick Katz and
Zeev Rudnick for stimulating discussions. 
The second-named author expresses his gratitude to the Institute for Advanced Study, where much of this work was written.


\section{Preliminaries}\label{sec:prelim}

We recommend the technical estimates in this section be omitted on a  first reading, and only referenced as needed in the proof, which begins in \secref{sec:setup}.

\subsection{Local Estimates}\

We begin with some elementary computations.
\begin{lem}\label{lem:sumsq}
For $p$ an odd prime,
$$
\#\{(x,y)\in\F_{p}^{2} \ : \ x^{2}+y^{2}=0\} \ = \
\twocase{}
{2p-1}{if $p\equiv1(4)$,}
{1,}{if $p\equiv3(4)$.}
$$
Moreover, for $\ell\neq0(p)$,
$$
\#\{(x,y)\in\F_{p}^{2} \ : \ x^{2}+y^{2}=\ell\} \ = \
\twocase{}
{p-1}{if $p\equiv1(4)$,}
{p+1,}{if $p\equiv3(4)$.}
$$
\end{lem}
\pf
Elementary.
\epf

For $\vp=\mattwos abcd,$ define
\be\label{eq:ffDef}
\ff(\vp) \ : = \ 
\tr (\vp {\, }^{\dag}\vp) \ = \
a^{2}+b^{2}+c^{2}+d^{2},
\ee
and for $\gep=\pm1$, set
\be\label{eq:rhoDef}
\rho(p) \ : = \ 
\frac1{|\SL_{2}(p)|}
\sum_{\g\in\SL_{2}(p)}
\bo_{\{\ff(\g)\equiv2\gep(p)\}}
.
\ee
Extend the definition of $\rho$ to all square-free $q$ by multiplicativity.
{\it A priori}, $\rho$ seems to depend on $\gep$, though the next lemma shows that it does not.
\begin{lem}\label{lem:rhoEval}
For $p$ an odd prime,
$$
\rho(p) \ = \ 
\twocase{}
{2p-1\over p(p+1)}{if $p\equiv1(4)$,}
{1\over p(p-1)}{if $p\equiv3(4)$.}
$$
Also, $\rho(2)=1/3$. 
\end{lem}
\pf
For $p=2$, two of the six matrices in $\SL_{2}(2)$ have $\ff=0$, so $\rho(2)=2/6$. Now assume $p\ge3$.
We need to count the number of $(a,b,c,d)\in\F_{p}^{4}$ with
$$
a^{2}+b^{2}+c^{2}+d^{2} \ = \ 2\gep,
\qquad\text{ and }\qquad
ad-bc \ = \ 1.
$$
We make the following linear change of variables:
\be\label{eq:varCh}
a \ = \ x+y,\quad
d \ = \ x-y,\quad
b \ = \ w+z,\quad
c \ = \ w-z,
\ee
which is invertible since $p\neq2$. The equations become
$$
x^{2}+y^{2}+z^{2}+w^{2} \ = \ \gep,
\qquad\text{ and }\qquad
x^{2}-y^{2}+z^{2}-w^{2} \ = \ 1,
$$
or equivalently,
\be\label{eq:varCh2}
x^{2}+z^{2} \ = \ 1+y^{2 }+ w^{2} \ = \
1+\bar 2(\gep-1)\ = \ 
\twocase{}
{1}{if $\gep=1$,}
{0}{if $\gep=-1$.}
\ee
Using \lemref{lem:sumsq} and $|\SL_{2}(p)|=p(p-1)(p+1)$ gives the claim.
\epf

Given $n\in\Z$, define $\Xi(q;n)$ on square-free $q$ by the expression
\be\label{eq:XiDef}
\Xi(p;n) \ : = \ \bo_{\{n\equiv0(p)\}} \ - \ \rho(p),
\ee
on primes $p$, and extend multiplicatively to $q$.

\begin{lem}\label{lem:XiSum0}
For any $\gw\in\SL_{2}(p)$ and $\gep=\pm1$,
$$
\frac1{|\SL_{2}(p)|}
\sum_{\g\in\SL_{2}(p)}
\Xi(p;\ff(\g\gw)-2\gep) \ = \ 0.
$$
\end{lem}
\pf
The coset $\gw$ plays no role since the $\g$ sum is over all of $\SL_{2}(p)$. The lemma follows from the definition \eqref{eq:rhoDef} of $\rho$.
\epf

The key estimate of this subsection is the following.

\begin{prop}\label{prop:key}
Let $\gw,\gw'\in\SL_{2}(p)$ and $\gep,\gep'\in\{\pm1\}.$ Then
\be\label{eq:key}
\frac1{|\SL_{2}(p)|}
\sum_{\g\in\SL_{2}(p)}
\Xi(p;\ff(\g\gw)-2\gep) 
\Xi(p;\ff(\g\gw')-2\gep') 
\ \ll \ 
\twocase{}
{\frac1p,}
{if $\gw\in\gw'\cdot\PO_{2}(p)$,}
{\frac1{p^{2}},}
{otherwise.}
\ee
Here we have defined
\be\label{eq:PO2Def}
\PO_{2}(p) 
\ : = \
\{ k\in\SL_{2}(p) \ : \
k {\,}^{\dag}k\equiv\pm I(p)\}
\ee
$$
\hskip1in
\  = \
\left\{\mattwo ab{-b}a \ : \
a^{2}+b^{2}\equiv\pm 1(p)\right\}.
$$
\end{prop}
\pf
Expanding $\Xi$ and using the definition \eqref{eq:rhoDef} of $\rho$, we must estimate
\be\label{eq:XiXi}
\frac1{|\SL_{2}(p)|}
\sum_{\g\in\SL_{2}(p)}
\bo_{\{\ff(\g\gw)\equiv2\gep(p)\}}
\bo_{\{\ff(\g\gw')\equiv2\gep'(p)\}}
\ - \ \rho(p)^{2}
\ee
The second term is plainly $\ll p^{-2}$ by \lemref{lem:rhoEval}.

If $p\equiv3(4)$, we may trivially bound $\bo_{\{\ff(\g\gw')\equiv2\gep'(p)\}}\le 1$, whence the first term is $\rho(p)=1/(p(p-1))\ll1/p^{2}$, as desired. Thus we may restrict to $p\equiv1(4)$.

If $\gw\in\gw'\cdot\PO_{2}(p)$, then $\ff(\g\gw)=\pm\ff(\g\gw')$, so if the signs $\gep,\gep'$ align, then the first term in \eqref{eq:XiXi} could be exactly $\rho(p)=(2p-1)/(p(p+1)) \asymp 1/p$. Thus we cannot do better than $1/p$ in this case. Now we seek extra cancellation when $\gw\not\in\gw'\cdot\PO_{2}(p)$.

Write 
$$
(\gw^{-1}\gw') {\, }^{\dag}(\gw^{-1}\gw')
\  =: \ 
\mattwo UVVW 
.
$$
Changing $\g\mapsto\g\gw^{-1}$ in \eqref{eq:XiXi} and using \eqref{eq:ffDef}, we must bound
$$
\frac1{|\SL_{2}(p)|}
\sum_{\g\in\SL_{2}(p)}
\bo_{\{\ff(\g)\equiv2\gep\}}
\bo_{\left\{\tr\left(\g{\,}^{\dag}\g \mattwos UVVW\right)\ \equiv\ 2\gep'\right\}}
.
$$
Writing $\g=\mattwos abcd$, the last equation becomes
$$
U(a^{2}+b^{2})+2V(ac+bd)+W(c^{2}+d^{2})\ \equiv\ 2 \gep'.
$$
Apply the same change of variables as in \eqref{eq:varCh}; then the equations become \eqref{eq:varCh2} and:
\be\label{eq:UVW}
U(1+2(xy+zw))+4V(xz-yw)+W(1-2(xy+zw)) \ \equiv \ 2\gep'.
\ee

Now suppose $\gep=1$ (the case $\gep=-1$ being similar). Then $x^{2}+z^{2}=1$ and $p\equiv1(4)$, so there are $p-1$ choices of $(x,z)$ by \lemref{lem:sumsq}. With $(x,z)$ fixed, \eqref{eq:UVW} becomes linear in $(y,w)$; we isolate $y$:
$$
2\big[(U-W)x-2Vw\big]y
\ \equiv \ 
2\gep'-(U+W)
-4Vxz
-2(U-W)zw
.
$$

Square and add $\bigg(2\big[(U-W)x-2Vw\big]w\bigg)^{2}$ to both sides  to take advantage of $y^{2}+w^{2}=0$. This gives a quartic equation in $w$ with everything else determined:
\beann
0&=&
\bigg[
2\gep'-(U+W)
-4Vxz
-2(U-W)zw
\bigg]^{2}
\\
&&
\hskip1.4in
+
\bigg(2\big[(U-W)x-2Vw\big]w\bigg)^{2}
\\
&=&
\bigg[
2\gep'-(U+W)
-4Vxz
\bigg]^{2}
-
4\bigg[
2\gep'-(U+W)
-4Vxz
\bigg]
(U-W)zw
\\
&&
+
4
(U-W)^{2}
w^{2}
-
16(U-W)Vxw^{3}
+16V^{2}w^{4}.
\eeann
This equation has at most 4 solutions in $w$, unless all the coefficients vanish, in which case $V=0$ and $U=W$. But $\det\mattwos UVVW=UW-V^{2}=1$, so $V=0$ implies $U=\overline W=W$. Hence $U=W=\pm1$, which means $\gw^{-1}\gw'\in\PO_{2}(p)$.
Since we have already dealt with this case, we may assume that
 the coefficients do not all vanish, whence there are at most 4 choices for $w$, from which $y$ is determined. In summary, there are $\ll p$ choices for $(x,z)$ and a bounded number of choices of $(y,w)$, while $|\SL_{2}(p)|\asymp p^{3}$; the ratio is $\ll 1/p^{2}$, as claimed.
\epf

\subsection{Spectral and Automorphic Estimates}\

We import here some lemmata from \cite{BourgainKontorovich2015}, first an automorphic estimate in $\SL_{2}(\Z)$.

\begin{lem}[{\cite[Lem. 2.13]{BourgainKontorovich2015}}]\label{lem:vfX}
Let $X\gg1$ be an increasing parameter. Then there is smooth bump function $\vf_{X}:\SL_{2}(\R)\to\R_{\ge0}$ with the following properties:

\begin{itemize}
\item It gives support to the norm-$X$ ball: If $\|g\|:=\sqrt{\tr(g{\,}^{\dag}g)}<X$, then
\be\label{eq:vfXge1}
\vf_{X}(g)\ \ge\ 1. 
\ee 

\item Furthermore,
\be\label{eq:vfGll1}
\sum_{\g\in\SL_{2}(\Z)}\vf_{X}(\g)  \ \ll \ X^{2}.
\ee

\item Finally, $\vf_{X}$ is evenly distributed in progressions: For any square-free $q$ and any $\g_{0}\in\SL_{2}(q)$,
\be\label{eq:vfModq}
\sum_{\g\in\SL_{2}(\Z)\atop \g\equiv\g_{0}(q)}
\vf_{X}(\g)
 \ = \
 \frac{1}{|\SL_{2}(q)|}
\sum_{\g\in\SL_{2}(\Z)}
\vf_{X}(\g)
\ + \
O(X^{3/2}).
\ee
\end{itemize}
All implied constants above are absolute. 
\end{lem}

\begin{rmk}
The error $X^{3/2}$ in \eqref{eq:vfModq} comes from using Selberg's $3/16$ spectral gap \cite{Selberg1965}; we are striving for simplest explicit exponents here, not optimal ones, so do not bother using best available exponents.
\end{rmk}

Finally, we will need SuperApproximation in our thin semigroup $\G_{\cA}$.
As discussed in \rmkref{rmk:AbsSpecGap}, we need this spectral gap to be absolute, so pick a fixed parameter $\cA_{0}=2$; then $\G_{\cA_{0}}=\G_{2}$.

\begin{lem}\label{lem:aleph}
For any $Y\gg1$, there is a non-empty subset 
$$
\aleph\ \subset \ \{\g\in\G_{2}:\|\g\|<Y\}
$$
and ``spectral gap''
\be\label{eq:gT}
\gT \ > \ 0,
\ee
so that, for any 
$q$ and any $\fa_{0}\in\SL_{2}(q)$,
\be\label{eq:alephEqui}
\#\{\fa\in\aleph \ : \ \fa\equiv\fa_{0}(q)\}
\
=
\
{1\over|\SL_{2}(q)|}
|\aleph|
+O(|\aleph|
{q^{C}Y^{-\gT}})
.
\ee
Here
$C$, $\gT$, and the implied constant are all absolute.
\end{lem}

\pf
A nearly identical statement is proved in \cite[Prop. 2.9]{BourgainKontorovich2015} with a weaker error term. The main ingredient there is a ``Prime Number Theorem''-type resonance-free region as proved in \cite{BourgainGamburdSarnak2011}. Now a resonance-free strip is  available (and does not require $q$ to be square-free) due to Magee-Oh-Winter/Bourgain-Kontorovich-Magee \cite{MageeOhWinter2016, BourgainKontorovichMagee2015}; substituting this result into the proof of \cite[Prop. 2.9]{BourgainKontorovich2015} gives the above claim.
\epf

\begin{rmk}
Actually the weaker statement \cite[Prop. 2.9]{BourgainKontorovich2015} using only \cite{BourgainGamburdSarnak2011} would already suffice for our purposes, see the treatment in \cite{BourgainKontorovich2015}. The resonance-free strip 
slightly simplifies the exposition, so we use it here.
\end{rmk}

\newpage

\section{%
Construction of $\Pi$ and 
the Sieving Theorem
}\label{sec:setup}

\subsection{Construction of the set $\Pi$}\

We  create here a certain subset $\Pi \ \subset \G_{\cA}$, all elements $\vp\in\Pi$ being of size $\|\vp\|\ll N$, a growing parameter, with $\Pi$ exhibiting a mutli-linear structure. First we break the parameter $ N$ as 
\be\label{eq:XYZ}
XYZ\ = \ N
,
\ee
and take the set $\aleph$ from \lemref{lem:aleph} with parameter $Y$.

The elements $\g\in\G_{\cA}$ of size $\|\g\|<X$ all have wordlength $\ell(\g)\asymp \log X$ in the generators \eqref{eq:GcADef}. By the Pigeonhole Principle, there is therefore a subset $\gW_{X}$ of $\G_{\cA}\cap B_{X}$ of size 
\be\label{eq:gWXsize}
\# \gW_{X} \ \gg \ { X^{2\gd}\over \log X},
\ee
(cf. \eqref{eq:Hensley})
all having the same wordlength. (We henceforth write $\gd$ for $\gd_{\cA}$, treating $\cA$ as fixed.) In the same way we construct the set $\gW_{Z}$ to parameter $Z$.

Then the set
\be\label{eq:PiDef}
\Pi\  := \ \gW_{X}\cdot\aleph\cdot \gW_{Z},
\ee
is a genuine subset (as opposed to multi-set) of $\G_{\cA}$, since each 
$$
\vp=\xi\cdot\fa\cdot\gw,
\qquad
\threecase
{\xi\in\gW_{X},}{}
{\fa\in\aleph,}{}
{\gw\in\gW_{Z}}{}
 $$ 
 is uniquely represented.

\subsection{The Sieving Theorem}\

In light of Lemmata \ref{lem:fundGeo} and \ref{lem:recipGeo}, 
we define $\Pi_{AP}$ to be the set of $\vp\in\Pi$ for which $\tr(\vp{}^{\dag}\vp)^{2}-4$ has no small prime factors, 
\be\label{eq:PiAPdef}
\Pi_{AP} \ := \
\{
\vp\in\Pi \ : \ 
p\mid (\tr(\vp^{\dag}\vp)^{2}-4)
\Longrightarrow
p>N^{1/350}
\}
.
\ee

An easy consequence of the main Sieving Theorem stated below is the following
\begin{thm}\label{thm:APs}
For any small  $\eta>0$, there is an $\cA=\cA(\eta)$, sufficiently large, and a choice of parameters $X,Y,Z$ in \eqref{eq:XYZ} so that
\be\label{eq:PiAPsize}
\#\Pi_{AP} \ > \ N^{2\gd-\eta},
\ee
as $N\to\infty$.
\end{thm}

The aforementioned Sieving Theorem is the following ``level of distribution'' result.

Recalling  $\ff$  defined in \eqref{eq:ffDef},
our sifting sequence is
$\fA=\{a_{N}(n)\}$ with
$$
a_{N}(n) \ := \
\sum_{\vp\in\Pi}\bo_{\ff(\vp)^{2}-4=n}.
$$
Note that $\fA$ is supported on $n<T$, where
\be\label{eq:TtoN}
T\asymp N^{4}. 
\ee
For square-free $\fq\ge1$, write
$$
|\fA_{\fq}| \ := \ 
\sum_{n\equiv0(\fq)}a_{N}(n)
,
$$
which measures the distribution of $a_{N}$ on multiplies of $\fq$.
\begin{thm}[The Sieving Theorem]\label{thm:Sieve}
For any small $\eta>0$, there is a sufficiently large $\cA=\cA(\eta)$ and a choice of the parameters $X,Y,Z$ so that 
the following holds.
Given a square-free $\fq$, there is a decomposition
\be\label{eq:cAfqDecomp}
|\fA_{\fq}|  \ = \ \gb(\fq)\cdot |\Pi| + r(\fq).
\ee
The function $\gb$ is multiplicative, and satisfies the ``quadratic sieve'' condition:
\be\label{eq:gbSieve}
\prod_{w\le p < z}\big(1-\gb(p)\big)^{-1}
\ \le \
C\cdot \left(
{\log z\over \log w}
\right)^{2}
.
\ee
Moreover, the ``remainder'' term $r(\fq)$ is controlled by:
\be\label{eq:rBnd}
\sum_{\fq<\cQ\atop\text{squarefree}}|r(\fq)|
\ \ll_{K} \ 
{|\Pi|\over \log^{K}N}
,
\qquad\text{ for any $K<\infty$,}
\ee
where the ``level of distribution'' $\cQ$ can be taken as large as
\be\label{eq:cQlevel}
\cQ \ = \ T^{1/72-\eta}.
\ee
Finally, the set $\Pi$ is large,
\be\label{eq:PiSize}
|\Pi| \ > \ N^{2\gd-\eta}.
\ee
\end{thm}

The deduction of 
\thmref{thm:APs} from \thmref{thm:Sieve} is completely standard, so we give a quick

\pf[Sketch]\
The sifting sequence $\fA$ has ``sieve dimension'' $\gk=2$, and any exponent of distribution $\ga<1/72$.
Taking $\ga=1/73$, say (again, we are not striving for optimal exponents), and using the crudest Brun sieve, see, e.g. \cite[Theorem 6.9]{FriedlanderIwaniecBook}, one shows that
\be\label{eq:aNsieve}
\sum_{n\atop (n,P_{z})=1}a_{N}(n) \ \gg \  {|\Pi|\over (\log N)^{2}},
\ee
where $P_{z}=\prod_{p<z}p$ and $z$ does not exceed $T^{\ga/(9\gk+1)}=T^{1/1387}=N^{4/1387}$; we take $z=N^{4/1400}=N^{1/350}$.
Of course any $n=\tr(\vp)^{2}-4$ coprime to $P_{z}$ has no prime factors below $z$. 
Then \eqref{eq:aNsieve} and \eqref{eq:PiSize}
confirm
 \eqref{eq:PiAPsize} after renaming constants.
\epf
We focus henceforth on establishing \thmref{thm:Sieve}. 

\subsection{The Decomposition and Dispersion}\

To prepare the proof, write,
for square-free $\fq\ge1$, 
$$
|\fA_{\fq}| \ := \ 
\sum_{n\equiv0(\fq)}a_{N}(n)
\ 
= \
\sum_{\tau\mod\fq\atop\tau^{2}\equiv4(\fq)}
\sum_{\vp\in\Pi}
\bo_{\{\ff(\vp)-\tau\equiv0(\fq)\}}
.
$$
To apply the ``dispersion'' method, we write
$$
\bo_{n\equiv0(p)} \ = \ \Xi(p;n) + \rho(p),
$$
with $\rho$ 
and $\Xi$ defined in \eqref{eq:rhoDef} and \eqref{eq:XiDef}, respectively.

Then
\bea\nonumber
|\fA_{\fq}|
 & = &
\sum_{\tau\mod\fq\atop\tau^{2}\equiv4(\fq)}
\sum_{\vp\in\Pi}
\prod_{p\mid\fq}
\bigg(
\Xi(p;\ff(\vp)-\tau)+\rho(p)
\bigg)
\\
\label{eq:fAfq}
& = &
\sum_{q\mid\fq}
\sum_{\tau\mod\fq\atop\tau^{2}\equiv4(\fq)}
\sum_{\vp\in\Pi}
\Xi(q;\ff(\vp)-\tau)\rho\left(\frac\fq q\right)
\eea

To give a decomposition towards \eqref{eq:cAfqDecomp}, we break the sum
\be\label{eq:fAfqBreak}
|\fA_{\fq}|  \ = \ \cM_{\fq}+r(\fq)
\ee
 according to whether $q<Q_{0}$ or not. 
 The two contributions are dealt with separately in the next two sections.


\section{Main Term Analysis}\label{sec:main}

From the decomposition \eqref{eq:fAfqBreak} of $\fA_{\fq}$ in  \eqref{eq:fAfq} the ``main'' term is:
\be\label{eq:cMis}
\cM_{\fq} \   = \
\sum_{q\mid \fq\atop q<Q_{0}}
\sum_{\tau(\fq)\atop\tau^{2}\equiv4}
\sum_{\vp\in\Pi}
\Xi(q;\ff(
\vp
)-\tau)
\rho\left(\frac\fq q\right)
.
\ee
The main goal of this section is to prove the following
\begin{thm}\label{thm:cMeval}
\be\label{eq:cMfqDecomp}
\cM_{\fq}  \ = \
\gb(\fq) \ |\Pi| +
r^{(1)}(\fq) ,
\ee
where $\gb$ is a multiplicative function defined on the primes by
\be\label{eq:gbIs}
\gb(p) \ : = \ 
\threecase
{\frac13,}{if $p=2$,}
{{2(2p-1)\over p(p+1)},}{if $p\equiv1(4)$,}
{{2\over p(p-1)},}{if $p\equiv3(4)$,}
\ee
and the ``remainder'' term $r^{(1)}$ satisfies:
\be\label{eq:r1Bnd}
\sum_{\fq<\cQ}
|r^{(1)}(\fq)| 
\ \ll \
|\Pi| \cQ^{\vep} 
Q_{0}^{C}Y^{-\gT}.
\ee
\end{thm}

To begin the proof, insert the construction \eqref{eq:PiDef} of $\Pi$ into \eqref{eq:cMis}, writing $\vp=\xi\fa\gw$. 
Since $\Xi(q;*)$ only depends on $*$ mod $q$, we decompose the $\fa$ sum along progressions mod $q$.
\beann
\cM_{\fq} 
&   = &
\sum_{q\mid \fq\atop q<Q_{0}}
\rho\left(\frac\fq q\right)
\sum_{\tau(\fq)\atop\tau^{2}\equiv4}
\sum_{\xi\in\gW_{X}}
\sum_{\gw\in\gW_{Z}}
\sum_{\fa_{0}\in\SL_{2}(q)}
\Xi(q;\ff(\xi\fa_{0}\gw)-\tau)
\left[
\sum_{\fa\in\aleph}
\bo_{\fa\equiv\fa_{0}(q)}
\right]
\eeann
 and apply expansion \eqref{eq:alephEqui}:
\beann
&   = &
\cM^{(1)}_{\fq}
+
r^{(1)}(\fq),
\eeann
where
\be\label{eq:cM1Is}
\cM^{(1)}_{\fq} \ := \
|\Pi|
\sum_{\tau(\fq)\atop\tau^{2}\equiv4}
\sum_{q\mid \fq\atop q<Q_{0}}
\rho\left(\frac\fq q\right)
\left[
\frac1{|\SL_{2}(q)|}
\sum_{\fa_{0}\in\SL_{2}(q)}
\Xi(q;\ff(\fa_{0})-\tau)
\right]
\ee
and
\beann
|r^{(1)}(\fq)| 
& \ll &
\sum_{q\mid \fq\atop q<Q_{0}}
\rho\left(\frac\fq q\right)
\sum_{\tau(\fq)\atop\tau^{2}\equiv4}
\sum_{\xi\in\gW_{X}}
\sum_{\gw\in\gW_{Z}}
\sum_{\fa_{0}\in\SL_{2}(q)}
|\Xi(q;\ff(\xi\fa_{0}\gw)-\tau)|
|\aleph |
q^{C}
Y^{-\gT}
\\
& \ll &
{\fq^{\vep}\over \fq}\
|\Pi|\
Q_{0}^{C}\
Y^{-\gT}
.
\eeann
Here we used $|\Xi|\le1$ and 
\lemref{lem:rhoEval} that $\rho(q)\ll q^{\vep}/q$.
Then \eqref{eq:r1Bnd} is immediately satisfied. 

Returning to $\cM^{(1)}_{\fq}$ in \eqref{eq:cM1Is}, the bracketed term vanishes unless $q=1$ by  \lemref{lem:XiSum0}, so 
we are left with
\beann
\cM^{(1)}_{\fq} \ = \
|\Pi|
2^{\nu(\fq)-\bo_{\{2\mid\fq\}}}
\rho(\fq)
.
\eeann
Here we elementarily evaluated the  contribution from the $\tau$ summation (see \cite[Lemma 4.1]{BourgainKontorovich2015}).
Inserting \lemref{lem:rhoEval}, we see that \eqref{eq:gbIs} is verified, 
completing the proof of
\thmref{thm:cMeval}.

%

\section{Error Term Analysis}\label{sec:err}

 The remainder term from \eqref{eq:fAfqBreak} is
$$
r(\fq) \ := \ 
\sum_{q\mid\fq\atop q>Q_{0}}
\sum_{t\mod q\atop t^{2}\equiv4(q)}
\sum_{\vp\in\Pi}
\Xi(q;\ff(\vp)-t)
\sum_{\tau\mod\fq\atop\tau^{2}\equiv4(\fq),\tau\equiv t (q)}
\rho\left(\frac\fq q\right)
$$
and total error is
$$
\cE \ := \
\sum_{\fq<\cQ}
|r(\fq)|,
$$
and we need to save a little more than $\cQ$ off of the trivial bound. 

The goal of this section is to prove the following
\begin{thm}\label{thm:cEbnd}
\be\label{eq:cEbnd}
\cE\ \ll\
T^{\vep}
|\Pi|
(XZ)^{1-\gd}
\left[
{\cQ^{4}\over X^{1/4}}
+
\frac1{Q_{0}^{1/2}}
+
{\cQ^{1/2}\over Z^{1/4}}
\right]
.
\ee
\end{thm}

First
write $\cE$ as
$$
\cE
\
=
\
\sum_{\fq<\cQ}
\gz(\fq)
r(\fq)
,
$$
where $\gz(\fq)=\overline{r(\fq)}/|r(\fq)|=\sgn r(\fq)$. 
Expanding gives
$$
\cE \ = \
\sum_{Q_{0}<q<\cQ}
\sum_{t\mod q\atop t^{2}\equiv4(q)}
\sum_{\vp\in\Pi}
\Xi(q;\ff(\vp)-t)
\gz_{1}(q,t)
$$
where
$$
\gz_{1}(q,t) \ := \
\sum_{\fq<\cQ\atop\fq\equiv0(q)}
\gz(\fq)
\sum_{\tau\mod\fq\atop\tau^{2}\equiv4(\fq),\tau\equiv t (q)}
\rho\left(\frac\fq q\right)
\ \ll
\
T^{\vep}
\sum_{\fq <\cQ/q}
{1\over \fq }
\ \ll
\
T^{\vep}.
$$
Decomposing $\Pi$ as $\gW_{X}\aleph\,\gW_{Z}$ gives
$$
\cE
\ = \
\sum_{Q_{0}<q<\cQ}
\sum_{t\mod q\atop t^{2}\equiv4(q)}
\sum_{\g\in\gW_{X}}
\sum_{\fa\in\aleph}
\sum_{\gw\in\gW_{Z}}
\Xi(q;\ff(\g\fa\gw)-t)
\gz_{1}(q,t),
$$
$$
\ \ll \
\sum_{Q_{0}<Q<\cQ\atop\text{dyadic}}
\sum_{\fa\in\aleph}
|\cE_{1}(\fa,Q)|
,
$$
where
$$
\cE_{1}(\fa,Q)
\ := \
\sum_{q\asymp Q}
\sum_{t\mod q\atop t^{2}\equiv4(q)}
\sum_{\g\in\gW_{X}}
\sum_{\gw\in\gW_{Z}}
\Xi(q;\ff(\g\fa\gw)-t)
\gz_{1}(q,t)
.
$$

\thmref{thm:cEbnd} follows immediately from 
\begin{prop}\label{prop:cE1bnd}
\be\label{eq:cE1bnd}
|\cE_{1}(\fa,Q)| \ \ll \
T^{\vep}
|\gW_{X}| 
|\gW_{Z}| 
(XZ)^{1-\gd}
\left[
{Q^{4}\over X^{1/4}}
+
\frac1{Q^{1/2}}
+
{Q^{1/2}\over Z^{1/4}}
\right]
.
\ee
\end{prop}

To begin the proof, apply the Cauchy-Schwarz inequality in the $\g$ variable and insert the smooth bump function $\vf_{X}$ from \lemref{lem:vfX}:
\bea\nonumber
|\cE_{1}(\fa,Q)|^{2} & \ll &
|\gW_{X}|\cdot
\sum_{\g\in\SL_{2}(\Z)}
\vf_{X}(\g)
\left|
\sum_{q\asymp Q}
\sum_{t\mod q\atop t^{2}\equiv4(q)}
\sum_{\gw\in\gW_{Z}}
\Xi(q;\ff(\g\fa\gw)-t)
\gz_{1}(q,t)
\right|^{2}
\\
\label{eq:cE12sum}
&\ll &
|\gW_{X}|\cdot
T^{\vep}
\sum_{q,q'\asymp Q}
\sum_{t\mod q\atop t^{2}\equiv4(q)}
\sum_{t'\mod q'\atop (t')^{2}\equiv4(q')}
\sum_{\gw,\gw'\in\gW_{Z}}
\\
\nonumber
&&
\times
\left|
\sum_{\g\in\SL_{2}(\Z)}
\vf_{X}(\g)
\Xi(q;\ff(\g\fa\gw)-t)
\Xi(q';\ff(\g\fa\gw')-t')
\right|
.
\eea
Having applied Cauchy-Schwarz,
we now need to save a little more than $Q^{2}$.
We first address the innermost $\g$ sum.

\begin{lem}\label{lem:gXiBnd}
Let
$$
\fq_{1} \ = \ 
\fq_{1}(\gw,\gw';q)
 \ : = \ 
\max_{\pm}(
\gcd(q, (\gw^{-1}\gw')\ {}^{\dag}(\gw^{-1}\gw') \mp I)),
$$
so that $\fq_{1}\mid q$ is the largest modulus for which $\gw^{-1}\gw'\in \PO_{2}(\fq_{1})$, the group defined
in \eqref{eq:PO2Def}.
Then
\bea\label{eq:gXiBnd}
&&
\left|
\sum_{\g\in\SL_{2}(\Z)}
\vf_{X}(\g)
\Xi(q;\ff(\g\fa\gw)-t)
\Xi(q';\ff(\g\fa\gw')-t')
\right|
\\
\nonumber
&& 
\hskip2in
\ll \
Q^{6}X^{3/2}
\ +\
\bo_{\{q=q'\}}
Q^{\vep}
X^{2}
{\fq_{1}\over q^{2}}
.
\eea
\end{lem}
\begin{rmk}
The first term above is a savings of $X^{1/2}$ against the loss of some powers of $Q$, which is more than the requisite $Q^{2}$ savings, as long as $Q$ is not too large relative to $X$. The second term is a savings of $Q$ from the $q=q'$ restriction, and a second factor of $Q^{2}/\fq_{1}$ from the $\fq_{1}/q^{2}$ term. If $\fq_{1}$ is small then this already saves more than $Q^{2}$, but if $\fq_{1}$ is of size $q$, then the net savings is $Q^{2}$ but no more. In that case, we will need just a bit extra savings from the fact that $\gw^{-1}\gw'\in\PO_{2}(\fq_{1})$ with such a large modulus $\fq_{1}$. 
\end{rmk}
\pf[Proof of \lemref{lem:gXiBnd}]
Let 
$$
\bar q:=[q,q']=\lcm(q,q'), \
\tilde q=(q,q')=\gcd(q,q'),\ 
q=q_{1}\tilde q,\  
q'=q_{1}'\tilde q,
$$ 
with $q_{1},$ $q_{1}'$ and $\tilde q$ pairwise coprime.
Because $\Xi(q,n)$ 
only depends on the residue of $n$ mod $q$, we break the innermost $\g$ sum into progressions, obtaining:
\bea\nonumber
\Bigg|
\sum_{\g\in\SL_{2}(\Z)}
\Bigg|
& = &
\sum_{\g_{0}\in\SL_{2}(\bar q)}
\Xi(q;\ff(\g_{0}\fa\gw)-t)
\Xi(q';\ff(\g_{0}\fa\gw')-t')
\left[
\sum_{\g\in\SL_{2}(\Z)\atop \g\equiv\g_{0}(\bar q)}
\vf_{X}(\g)
\right]
\\
\nonumber
& \ll &
X^{2}
\left|
\frac1{|\SL_{2}(\bar q)|}
\sum_{\g_{0}\in\SL_{2}(\bar q)}
\Xi(q;\ff(\g_{0}\fa\gw)-t)
\Xi(q';\ff(\g_{0}\fa\gw')-t')
\right|
\\ \label{eq:g0sum}
\\ \nonumber
&&
\hskip3in
+
O(\bar q^{3}X^{3/2})
,
\eea
where we used \eqref{eq:vfModq} and \eqref{eq:vfGll1}. 
Since  $\bar q\ll Q^{2}$, the last term contributes $Q^{6}X^{3/2}$ to \eqref{eq:gXiBnd}.

Now, the 
remaining
 $\g_{0}$ sum in \eqref{eq:g0sum} is multiplicative, so decomposing $\bar q=q_{1}q_{1}'\tilde q$, we can write it as:
\beann
\Bigg|
\sum_{\g_{0}\in\SL_{2}(\bar q)}
\Bigg|
&=&
\left|
\frac1{|\SL_{2}(q_{1})|}
\sum_{\g_{0}\in\SL_{2}(q_{1})}
\Xi(q_{1};\ff(\g_{0}\fa\gw)-t)
\right|
\\
&&\times
\left|
\frac1{|\SL_{2}( q_{1}')|}
\sum_{\g_{0}\in\SL_{2}( q_{1}')}
\Xi(q_{1}';\ff(\g_{0}\fa\gw')-t')
\right|
\\
&&\times
\left|
\frac1{|\SL_{2}(\tilde q)|}
\sum_{\g_{0}\in\SL_{2}(\tilde q)}
\Xi(\tilde q;\ff(\g_{0}\fa\gw)-t)
\Xi(\tilde q;\ff(\g_{0}\fa\gw')-t')
\right|
.
\eeann
From \lemref{lem:XiSum0}, we see that the first two terms completely vanish unless $q_{1}=q_{1}'=1$, that is, $q=q'=\tilde q=\bar q$.
For the third term, we apply the key \propref{prop:key}; then every $p\mid q$, contributes either $1/p$ or $1/p^{2}$, depending on whether $\gw^{-1}\gw'\in\PO_{2}(p)$ or not. 
This savings is exactly captured by $Q^{\vep}\fq_{1}/q^{2}$, completing the proof.
\epf

\pf[Proof of \propref{prop:cE1bnd}]
Inserting \eqref{eq:gXiBnd} into \eqref{eq:cE12sum} gives:
\bea\nonumber
|\cE_{1}(\fa,Q)|^{2} & \ll &
T^{\vep}
|\gW_{X}|
Q^{2}
|\gW_{Z}|^{2}
Q^{6}X^{3/2}
\\
\nonumber
&&
+
T^{\vep}
|\gW_{X}|
\sum_{q
\asymp Q}
\sum_{\gw
\in\gW_{Z}}
\sum_{\fq_{1}\mid q}
X^{2}
{\fq_{1}\over q^{2}}
\left[
\sum_{\gw'\in\SL_{2}(\Z)\atop \gw^{-1}\gw'\in\PO_{2}(\fq_{1})}\vf_{Z}(\gw')
\right]
,
\\
\label{eq:gwpSum}
\eea
where we extended the $\gw'$ sum to all of $\SL_{2}(\Z)$ and again inserted the bump function $\vf$ from \lemref{lem:vfX}. Now break
the innermost $\gw'$ sum in \eqref{eq:gwpSum} into progressions mod $\fq_{1}$, and apply \eqref{eq:vfModq} and \eqref{eq:vfGll1} to obtain:
\beann
\left[
\sum_{\gw'\in\SL_{2}(\Z)}
\right]
&=&
\sum_{\gw'_{0}\in\gw\cdot\PO_{2}(\fq_{1})}
\sum_{\gw'\in\SL_{2}(\Z)\atop \gw'\equiv\gw'_{0}(\fq_{1})}
\vf_{Z}(\gw')
\ \ll\ 
|\PO_{2}(\fq_{1})|
\left[
{Z^{2}\over \fq_{1}^{3}}
+
Z^{3/2}
\right]
\\
&\ll&
\fq_{1}^{\vep}
\left[
{Z^{2}\over \fq_{1}^{2}}
+
\fq_{1}Z^{3/2}
\right]
,
\eeann
since $|\PO_{2}(\fq_{1})|\ll \fq_{1}^{1+\vep}$.
The contribution 
of this
to \eqref{eq:gwpSum} is then
\beann
&\ll&
T^{\vep}
|\gW_{X}|
\sum_{q
\asymp Q}
\sum_{\gw
\in\gW_{Z}}
\sum_{\fq_{1}\mid q}
X^{2}
{\fq_{1}\over q^{2}}
\left[
{Z^{2}\over \fq_{1}^{2}}
+
\fq_{1}Z^{3/2}
\right]
\\
&\ll&
T^{\vep}
|\gW_{X}|
X^{2}
|\gW_{Z}|
Z^{2}
\left[
{1\over Q}
+
{Q\over Z^{1/2}}
\right]
.
\eeann
Combined with the first term of \eqref{eq:gwpSum} and  \eqref{eq:gWXsize}, this gives 
\eqref{eq:cE1bnd}, as claimed.
\thmref{thm:cEbnd} follows immediately.
\epf


\section{Proof of the Sieving Theorem 
}\label{sec:pfThm1}

We proceed now to prove \thmref{thm:Sieve}.
Combining 
\eqref{eq:cMfqDecomp} with \eqref{eq:fAfqBreak} gives the decomposition \eqref{eq:cAfqDecomp}. 
The content of \eqref{eq:gbSieve} is, roughly, that $\gb(p)\sim 2/p$ on average; indeed, from \eqref{eq:gbIs} we have that
$$
\gb(p) \ = \ \twocase{}
{\frac4p+O(p^{-2}),}{if $p\equiv1(4)$,}
{O(p^{-2}),}{if $p\equiv3(4)$,}
$$
so 
\eqref{eq:gbSieve}
is elementarily verified. 
Combining \eqref{eq:r1Bnd} and 
\eqref{eq:cEbnd} gives \eqref{eq:rBnd}, as long as the following inequalities are satisfied:
\bea\label{eq:y}
C\ga_{0} & <& \gT y,
\\
\label{eq:x}
4\ga+(1-\gd)(x+z) & <& x/4,
\\
\label{eq:ga0}
(1-\gd)(x+z) & <& \ga_{0}/2,
\\
\label{eq:z}
\ga/2+(1-\gd)(x+z) & <& z/4.
\eea
Here
$$
Q_{0}\ = \ N^{\ga_{0}},\
\cQ  \ = \ N^{\ga},\
X \ = \ N^{x},\
Y \ = \ N^{y},\
Z \ = \ N^{z}.
$$

\begin{rmk}
Treating $1-\gd$ as $0$ and $x+z$ as $1$, one quickly sees that the best one can do is the choice $\ga\approx 1/18$, $x\approx 8/9$, $z\approx 1/9$.
\end{rmk}

Let $\eta>0$ be given, and set 
$$
\ga \ = \ 1/18-\eta.
$$
Since $\cQ=N^{\ga}$ and $N\asymp T^{1/4}$ (see \eqref{eq:TtoN}), this gives the exponent of distribution $1/72$ claimed in \eqref{eq:cQlevel}.
Next we set
$$
x\ = \ \frac 89-\eta,
$$
and assume at first that $1-\gd<\eta$ (more stringent restrictions on $\gd$ will follow). Then since $x+z<1$, we have
$$
x\ = \ 16\ga+15\eta \   > \ 16 \ga + 4 \eta\  > \ 16\ga +4(1-\gd)(x+z).
$$
That is, \eqref{eq:x} is satisfied. Similarly, we set
$$
z\ = \ \frac 19-\eta,
$$
whence \eqref{eq:z} holds once $1-\gd<\eta/4$. This means that $y=2\eta$, so \eqref{eq:y} is satisfied when
$$
\ga_{0} \ = \ {\gT\eta\over C}.
$$
Finally, for \eqref{eq:ga0} to hold, we need
$$
\gd \ > \ 1- {\gT\eta\over 2C}\left(1-2\eta\right)^{-1}.
$$
Recalling that $\gd=\gd_{\cA}$, this can be achieved (cf. \eqref{eq:gdcA}) by taking $\cA$ sufficiently large.

\begin{rmk}\label{rmk:NeedAbsGap}
It is here that we are crucially using that the parameters  $\gT$ and $C$ coming from the ``spectral gap'' estimate \eqref{eq:alephEqui} are independent of $\cA$, and only depend on the fixed quantity $\cA_{0}=2$; that is, they are absolute.
\end{rmk}


\section{Proof of \thmref{thm:1}}\label{sec:pfMain}


\subsection{Proof of \thmref{thm:main}}\label{sec:pftr2}\

\begin{lem}
As $t\to\infty$,
\be\label{eq:mulBnd}
\#\{
\g\in\SL_{2}(\Z) 
 \ : \
\tr(\g^{\dag}\g)=t
\}
\ \ll \ 
t^{\vep}.
\ee
\end{lem}
\pf
One must count the number of $( a,b,c,d)\in\Z^{4}$ having 
$ad-bc=1$
and 
$a^{2}+b^{2}+c^{2}+d^{2}=t$.
 Changing variables to $a=x+y$, $d=x-y$, $b=w+z$, $c=w-z$ gives the equations
$x^{2}-y^{2}+z^{2}-w^{2}=1$, and $2(x^{2}+y^{2}+z^{2}+w^{2})=t$. That is, any solution to the former equations in integers gives one to the latter equations. It is elementary to see there are at most $t^{\vep}$ solutions to the latter.
\epf

In this section, write $\ga=1/350$, so that we can write \eqref{eq:PiAPdef} as
$$
\Pi_{AP}
\ = \
\{
\g\in\Pi\ : \ p\mid (\tr(\g^{\dag}\g)^{2}-4)
\ \Longrightarrow\  p>
N^{\ga}\}
.
$$
\thmref{thm:main} asks us to count
\bea
\label{eq:mainThm}
&&
\hskip-.5in
\#\{
\g\in\G_{\cA}\cap B_{N}\ : \ \tr(\g^{\dag}\g)^{2}-4\text{ is square-free}\}
\\
\nonumber
&
\ge
&
\#\{
\g\in\Pi_{AP} \ : \ \tr(\g^{\dag}\g)^{2}-4\text{ is square-free}\}
\\
\label{eq:nsqf}
&
>
&
N^{2\gd-\eta}
\ - \
\#
\Pi_{AP}^{\square}
,
\eea
where
we used \eqref{eq:PiAPsize} and defined
$$
\Pi_{AP}^{\square}
\ := \
\{
\g\in\Pi_{AP} \ : \ \tr(\g^{\dag}\g)^{2}-4\text{ is not square-free}\}
.
$$
Now, 
for each $\g\in\Pi_{AP}^{\square},$ there is a prime $p$ with 
$p^{2}\mid (\tr(\g^{\dag}\g)^{2}-4)$. Since $\g\in\Pi_{AP}$, we thus have that $p>N^{\ga}$, and
moreover, 
$p^{2}$ divides either $\tr(\g^{\dag}\g)+2$ or $\tr(\g^{\dag}\g)-2$; in particular, $p\ll N$. 
Therefore, reversing orders and applying \eqref{eq:mulBnd}, we have
\beann
\#
\Pi_{AP}^{\square}
&
\le
&
\sum_{N^{\ga}<p\ll N}
\sum_{t<N^{2}\atop t^{2}-4\equiv0(p^{2})}
\#\{\g\in\G_{\cA}\cap B_{N} \ : \ \tr(\g^{\dag}\g)=t\}
\\
&
\ll
&
\sum_{N^{\ga}<p\ll N}
{N^{2}\over p^{2}}
N^{\vep}
\
\ll
\
N^{2-\ga+\vep}
.
\eeann

Since $\ga=1/350$ is fixed, it is clear that by making $\gd=\gd_{\cA}$ sufficiently near $1$ (by taking $\cA$ large), one gets the desired main term from \eqref{eq:nsqf}. This completes the proof of \thmref{thm:main}.

\subsection{Proof of \thmref{thm:1}}\

Again, this will be an easy consequence of \thmref{thm:main}. 
Each $\g\in\G_{\cA}\cap B_{N}$ arising in \eqref{eq:mainThm} gives a hyperbolic conjugacy class $[\g^{\dag}\g]$ of trace at most $N^{2}$, for which the corresponding geodesic is low-lying (relative to $\cA$), reciprocal, and fundamental. The only two issues  are {\it (a)} that the class need not be primitive, and {\it (b)} different $\g$'s can give rise to the same geodesic. 
Since the wordlength metric is commensurate with the logarithm of the archimedean metric, the number of imprimitive classes (that is, $\g$ which, as symbols in the generators of $\G_{\cA}$, have a repeating sequence) is easily bounded by
$N^{1+\vep}$; these can safely be discarded from \eqref{eq:mainThm} without affecting the cardinality.
The latter {\it (b)} happens when the symbols generating $\g$ and $\g'$, say, are the same up to a cyclic permutation. This adds at most $\log N$ to the multiplicity of \eqref{eq:mainThm}, and can thus also be safely discarded.
In summary, we have produced $N^{2-\eta}$ low-lying, fundamental, and reciprocal closed geodesics with trace bounded by $N^{2}$, as claimed. This completes the proof.

\newpage

\bibliographystyle{alpha}

\bibliography{../../AKbibliog}

\end{document}